\documentclass[12pt]{amsart}

\usepackage{amssymb}
\usepackage{isolatin1}    
\usepackage{mathrsfs}

\begin{document}
\newtheorem{theorem}{Theorem}
\newtheorem{lemma}{Lemma}

\title{The support theorem for the complex Radon transform of distributions}

\author[Sekerin]{A.B. Sekerin}
\address{ The Orel State University, Komsomolskaya 95, Orel, 302026, Russia }
\email{sekerin@orel.ru}


\thanks{The author gratefully acknowledges support of the Deutscher Akademischer Austauschdienst (DAAD) for the visit at Mathematisches Institut der Heinrich-Heine-Universit\"at D\"usseldorf, where this research was carried out.} 

\begin{abstract}The complex Radon transform $\hat F$ of a rapidly decreasing distribution $F\in\mathscr{O}_C^{\prime}(\mathbb{C}^n)$ is considered.  
A compact set $K\subset\mathbb{C}^n$ is called linearly convex if the set $ \mathbb{C}^n \setminus K$ is a union of complex hyperplanes. Let  $\hat K$ denote the set of complex hyperplanes which meet $K$. The main result of the paper establishes the conditions on a linearly convex compact $K$ under which the support theorem for the  complex Radon transform is true: from the relation $\hbox{supp}(\hat F)\subset\hat K$  it follows that $F\in\mathscr{O}^{\prime}_C(\mathbb{C}^n)$ is compactly supported and $\hbox{supp}(F)\subset K$.      
\end{abstract}

\maketitle

If $f$ is the function defined on $\mathbb{R}^n$ ($\mathbb{C}^n$), the classical real (complex)  Radon transform $Rf$ of $f$ is the function defined on hyperplanes; the value of $Rf$ at a given hyperplane is the integral of $f$ over that hyperplane. For the theory of the Radon transform we refer to J. Radon \cite{Radon}, F. John \cite{John1}, \cite{John2}, I.M.Gel'fand, M.I.Graev, and N.Ya. Vilenkin \cite{Gelfand},
S. Helgason \cite{Helgason1},  \cite{Helgason2}, D. Ludwig \cite{Ludwig}, A. Hertle \cite{Hertle1}. One of the basic results on the classical Radon transform is Helgason's support theorem \cite{Helgason1}: A rapidly decreasing function must vanish outside a ball if its real  Radon transform does. This theorem holds for every convex compact set in $\mathbb{R}^n$ and remains valid for rapidly decreasing distributions \cite{Hertle1}.  

In the present paper we prove the support theorem for  the complex Radon transform of  distributions.

{\bf Notations}. \ For $z,w\in\mathbb{C}^n$ we write  $\langle z,w\rangle=\sum z_{j}w_{j}$.  $B^n(z,R):=\{w\in\mathbb{C}^n \left|\right. |w-z|<R\}$ denotes the euclidean ball of center $z$ and radius r in $\mathbb C^n$. If X is a set, we denote by $\bar X$ the closure of X. The standard Lebesgue measure in $\mathbb{C}^n$ is $d\omega_{2n}$. $S^{2n-1}$ denotes the unit sphere in $\mathbb{C}^n$, and $d\sigma$ is the area element on $S^{2n-1}$. For $n$-tuples $p=(p_1,p_2,\ldots,p_n)$ and $q=(q_1,q_2,\ldots,q_n)$ of non-negative integers, we denote by $\partial^p\bar\partial^q$ the partial derivative 
$$\frac{\partial^{|p|+|q|}}
{\partial z_1^{p_1}\ldots\partial z_n^{p_n}
\partial \bar z_1^{q_1}\ldots\partial\bar z_n^{q_n}}$$
of order $|p|+|q|=p_1+\ldots+p_n+ q_1+\ldots+q_n$.     
Similarly, for $z=(z_1,\ldots,z_n)$ we write $z^p=z_1^{p_1}\ldots  z_n^{p_n}$, 
$\bar z^q=\bar z_1^{p_1}\ldots  \bar z_n^{q_n}$.  
For a domain $\Omega\subset\mathbb{C}^n$, we denote by    
$\mathscr{S}(\Omega)$, $\mathscr{D}(\Omega)$, and $\mathscr{E}(\Omega)$ the spaces of rapidly decreasing $C^{\infty}$ functions, $C^{\infty}$ functions with compact support, and $C^{\infty}$ functions, respectively. The dual spaces $\mathscr{S}^{\prime}(\Omega)$, $\mathscr{D}^{\prime}(\Omega)$,  and $\mathscr{E}^{\prime}(\Omega)$ are the spaces of tempered distributions, distributions,  and distributions with compact support, respectively.  

If $\varphi\in\mathscr{S}(\mathbb{C}^n)$, the standard complex Radon transform of $\varphi$ (denoted by $\hat\varphi$) is defined by 
\begin{equation}\label{g2eq1}
\hat\varphi(\xi,s)=\frac{1}{|\xi|^2}
\int\limits_{\langle z,\xi\rangle=s}\varphi(z)\,d\lambda(z),
\end{equation}
where $(\xi,s)\in(\mathbb{C}^{n}\setminus 0)\times\mathbb{C}$, and  $d\lambda(z)$ is the area element on the hyperplane $\{z: \langle z,\xi\rangle=s\}$. For a set $A\subset\mathbb{C}^n$, we denote by $\hat A$ the set of all $(\xi,s)\in(\mathbb{C}^{n}\setminus 0)\times\mathbb{C}$ such that the hyperplane $\{z: \langle z,\xi\rangle=s\}$ meets $A$. A set $A\subset\mathbb{C}^n$ is called linearly convex if, for every $w\notin A$, there is a complex hyperplane $\{z: 
\langle z,\xi\rangle=s\}$ which contains $w$ and does not meet $A$ (see Martineau \cite{Martineau}).

We use the approach of Gel'fand et al. \cite{Gelfand} to introduce the complex Radon transform of distributions. 
Let $X=S^{2n-1}\times\mathbb{C}$, and let $\mathscr{E}(X)$ be the set of complex-valued functions $\varphi(w,s)$ on $S^{2n-1}\times\mathbb{C}$ which satisfy the following conditions:
\begin{itemize}
\item[(a)]
Functions  $\varphi(w,s)$ are infinitely differentiable with respect to $s$.  
\item[(b)]
For all $p,q\ge 0$ the derivatives 
$$\frac{\partial^{p+q}\varphi(w,s)}{\partial s^p\partial\bar s^q}$$
are continuous on $S^{2n-1}\times\mathbb{C}$. 
\item[(c)]
$\varphi(we^{i\theta},se^{i\theta})=\varphi(w,s)$ for all $\theta\in[0,2\pi]$. 
\end{itemize}
We give $\mathscr{E}(X)$ the topology defined by the system of seminorms
$$q_k(f)=\max\limits_{k_1+k_2\le k}\max\limits_{|s|\le k}\max
\limits_{w\in S^{2n-1}}
\left|\frac{\partial^{k_1+k_2}f(w,s)}{\partial s^{k_1}\partial\bar s^{k_2}}\right|.$$
By $\mathscr{D}(X)$ we denote  the space of all compactly supported functions in $\mathscr{E}(X)$. We give $\mathscr{D}(X)$ the standard topology of the inductive limit of the spaces 
$$\mathscr{D}_m=\left\{\varphi\in\mathscr{E}(X): \hbox{supp}(\varphi) \subset
S^{2n-1}\times \{|s|\le m\}\right\}.$$
Let $R\mathscr{D}(X)$ be the subspace of $\mathscr{D}(X)$ formed by the Radon transforms $\hat\varphi$ of functions in $\mathscr{D}(\mathbb{C}^n)$ (the equality $\hat\varphi(we^{i\theta},se^{i\theta})\equiv\hat\varphi(we^{i\theta},se^{i\theta})$ follows for 
$\varphi\in\mathscr{D}(\mathbb{C}^n)$ from the definition of $\hat\varphi$). Similarly, 
we define the subspace $R\mathscr{S}(X)$ of $\mathscr{S}(X)$. 

The dual Radon transform  is the operator $R^*: \mathscr{E}(X)\rightarrow\mathscr{E}(\mathbb{C}^n)$   given by 
$$
[R^*(f)](z)=\int\limits_{S^{2n-1}}f(w,\langle z, w\rangle)\,d\sigma(w).$$
It is easy to see that the operator $R^*$ is continuous. It follows from the definition of the Radon transform that 
\begin{equation}\label{g2eq2}
\int\limits_{\mathbb{C}^n}[R^*(f)](z)\varphi(z)\, d\omega_{2n}(z)=
\int\limits_{\mathbb{C}}\int\limits_{S^{2n-1}}f(w,s)\hat\varphi(w,s)\,d\sigma(w)d\omega_2(s)
\end{equation}
for every function $\varphi\in\mathscr{D}(\mathbb{C}^n)$. 

Let $M_{\mathscr{D}}$ be the subspace of $\mathscr{D}(X)$ formed by the functions 
\begin{equation}\label{g2eq2a}
\psi(w,s)=\frac{\partial^{2n-2}\hat\varphi(w,s)}{\partial s^{n-1}\partial\bar s^{n-1}}, \quad \hat\varphi\in R\mathscr{D}(X).
\end{equation}
We give $M_{\mathscr{D}}$ the topology induced from $\mathscr{D}(X)$. 

{\bf Definition 1.} {\it Let $F\in\mathscr{D}^{\prime}$. The Radon transform $RF$ of $F$ is the functional on $M_{\mathscr{D}}$ given by}  
\begin{equation}\label{g2eq4}
\langle RF,\psi\rangle = \langle F,R^*\psi\rangle.
\end{equation}

For every function $\varphi\in\mathscr{S}(\mathbb{C}^n)$ the following inversion formula holds \cite[p. 118]{Gelfand}: 
\begin{equation}\label{g2eq3}
\varphi(z)=(-1)^{n-1}c_nR^*\left(\frac{\partial^{2n-2}\hat\varphi(w,s)}{\partial s^{n-1}\partial\bar s^{n-1}}\right), 
\end{equation}
where $\hat\varphi(w,s)$ is the Radon transform of $\varphi$, and $c_n>0$. 
It follows from the inversion formula (\ref{g2eq3}) that for each function $\psi\in M_{\mathscr{D}}$ the function $R^*(\psi)(z)$ belongs to $\mathscr{D}(\mathbb{C}^n)$. 
Therefore the functional $RF$ is well defined.  

{\bf Definition 2.} {\it We say that the Radon transform $RF$ of a distribution $F\in\mathscr{D}^{\prime}$ is defined as a distribution if the functional $RF$ given by (\ref{g2eq4}) can be extended to a continuous functional on $\mathscr{D}(X)$}. 

It has been shown in \cite{Hertle1} that there are distributions in $\mathbb{R}^m$ whose real Radon transforms are not defined as distributions.  It is natural to suppose that there are such examples in the case of the complex Radon transform.  If the distribution $F$ is given by the function $f(z)\in\mathscr{S}(\mathbb{C}^n)$, then it follows from (\ref{g2eq3}) and (\ref{g2eq2}) that the Radon transform $RF$ is defined as a distribution and it is given by the function $\hat f(w,s)$.
  
We denote by $\mathscr{O}^{\prime}_C(\mathbb C^n)$ the space of rapidly decreasing distributions \cite[p. 419]{Horvath}. A distribution $T\in\mathscr{D}^{\prime}
(\mathbb{ C}^n)$ belongs to $\mathscr{O}^{\prime}_C(\mathbb{C}^n)$ if and only if 
for every $k\in\mathbb{Z}$ the distribution $(1+|x|^2)^kT$ is integrable; i.e.,
\begin{equation}\label{g2eq4a}
(1+|x|^2)^kT=\sum\limits_{|p|+|q|\le m(k)}\partial^p\bar\partial^q\mu_{pq}(k),
\end{equation}
where $m(k)\in\mathbb{N}$ and $\{\mu_{pq}\}(k)$ is a finite family of bounded measures on $\mathbb{C}^n$. In particular, every distribution with compact support is rapidly decreasing. 

Let $T\in\mathscr{O}^{\prime}_C(\mathbb C^n)$. We show that equality (\ref{g2eq4}) defines  the extension of the Radon transform $RT$ to a continuous linear functional on $\mathscr{D}(X)$.  Let $h(w,s)\in\mathscr{D}(X)$ be such that $|h(w,s)|\le 1$. There is $R>0$ such that $h(w,s)=0$ for $|s|\ge R$, and we have 

\begin{equation}\label{g2eq5}
\left|[R^*h](z)\right|\le\int\limits_{S^{2n-1}}|h(w,\langle z, w\rangle)|d\sigma(w)\le\int\limits_{|\langle z, w\rangle|\le R}d\sigma(w)\le d_n\max\left(1, \frac{R^2}{|z|^2}\right),
\end{equation}
where $d_n>0$. Suppose that the sequence $\{h_N(w,s)\}$ in $\mathscr{D}(X)$ converges to $0$. Then, for every multi-indices $p$ and $q$, we have

\begin{equation}\label{g2eq5b}
\partial^p\bar\partial^q\left[R^*(h_N)\right](z)=\int\limits_{S^{2n-1}}
\frac{\partial^{|p|+|q|}}{\partial s^{|p|}\partial\bar s^{|q|}}
h_N(w,\langle z,w \rangle)
w^p\bar w^qd\sigma(w).
\end{equation}
There exists $R>0$ such that $\hbox{supp}(h_N)\subset S^{2n-1}\times\{s: |s|\le R\}$ for  all $N$. Then it follows from (\ref{g2eq5}) and (\ref{g2eq5b}) that 
\begin{equation}\label{g2eq5a}
\left|\partial^p\bar\partial^q\left[R^*(h_N)\right](z)\right|\le
d_n\max\left(1, \frac{R^2}{|z|^2}\right)\max\limits_{w,s}\left|\frac{\partial^{|p|+|q|}}{\partial s^{|p|}\partial\bar s^{|q|}}
h_N(w,s)\right|.
\end{equation}
This means that the functions  $\left[R^*(h_N)\right](z)$, together with derivatives of all orders, vanish at infinity.  
By the definition of the topology of $\mathscr{D}(X)$ we have

\begin{equation}\label{g2eq6}
\lim\limits_{N\to\infty}\max\limits_{w,s}\left|\frac{\partial^{|p|+|q|}}{\partial s^{|p|}\partial\bar s^{|q|}}
h_N(w,s)\right|=0.
\end{equation}
We set $k=0$ in (\ref{g2eq4a}). Then we obtain from (\ref{g2eq4a}) and 
(\ref{g2eq4}) that  
$$\langle RT,h_N\rangle=\langle T,[R^*h_N]\rangle=
\sum\limits_{|p|+|q|\le m}(-1)^{|p|+|q|}\int\limits_{\mathbb{C}^n}
\partial^p\bar\partial^q[R^*h_N](z)d\mu_{pq}(z).$$
Since the measures $\mu_{pq}$ are bounded,  it follows from (\ref{g2eq5a}) and (\ref{g2eq6}) that $\langle RT,h_N\rangle\to 0$ as $N\to\infty$. Thus, for every 
$T\in\mathscr{O}^{\prime}_C(\mathbb{C}^n)$, the functional $RT$ is well-defined and continuous on $\mathscr{D}(X)$. 

\begin{theorem}\label{th1} Let $T\in\mathscr{O}_C^{\prime}(\mathbb{C}^n)$ and let $K\subset\mathbb{C}^n$ be a linearly convex compact set. Suppose
that for every $z\notin K$ there exists a hyperplane $P=\{\lambda:
\langle\lambda,w_0\rangle=s_0\}$ satisfying the following conditions:\\
(i) $P$ contains $z$.\\
(ii) $P$ does not meet $K$.\\
(iii) The set $\mathbb{C}\setminus K_{w_0}$ is connected, where
$K_{w_0}=\{\langle \lambda,w_0\rangle\}_{\lambda\in K}$ is the projection
of $K$ on $w_0$.
Then $T$ has support in $K$ if and only if its Radon transform $RT$ has support in $\hat K$. 
\end{theorem}

{\it Remark.}\ Theorem \ref{th1} was proved by the author in the special case in which the distribution $T$ is given by a compactly supported continuous function \cite{Sekerin}. The proof of Theorem \ref{th1} is based on the properties of the convolution of $T$ and smooth compactly supported functions. As in the proof of the  similar theorem for the real Radon transform and convex compact sets \cite{Hertle1}, the proof of Theorem \ref{th1} can be easily reduced to the case of regular distributions if for small enough $\varepsilon>0$ the set 
$$K_{\varepsilon}=\bigcup\limits_{z\in K}\bar B^n(z,\varepsilon)$$
also satisfies the conditions (i)-(iii). It should be noted that, in contrast to the case of convex compacts, there are examples of compact sets $K$ satisfying (i)-(iii) such that the set $K_{\varepsilon}$ does not satisfy the condition (iii) for every $\varepsilon>0$. Since it has been shown in \cite{Sekerin} that assumption (iii) in Theorem \ref{th1} is essential,  Theorem \ref{th1} is  not  a simple consequence of the result of \cite{Sekerin}.

{\it Proof of Theorem \ref{th1}.}\ Suppose that $T\in\mathscr{O}^{\prime}_C(\mathbb{C}^n)$ has support in $K$. Then $T\in\mathscr{E}(\mathbb{C}^n)$. 
Let $h(w,s)\in\mathscr{D}(X)$ be such that $\hbox{supp}(h)\subset X\setminus \hat K$. If $z\in K$, then  the point $(w,\langle z,w\rangle)$ belongs to $\hat K$ for every $w\in S^{2n-1}$. Therefore the functions 
$$[R^*h](z)=\int\limits_{S^{2n-1}}h(w,\langle z, w\rangle) d\sigma(w), $$
$$\partial^p\bar\partial^q\left[R^*(h)\right](z)=\int\limits_{S^{2n-1}}
\frac{\partial^{|p|+|q|}}{\partial s^{|p|}\partial\bar s^{|q|}}
h(w,\langle z,w \rangle)
w^p\bar w^qd\sigma(w)
$$
vanish on $K$. So $[R^*h](z)$ is an infinitely differentiable function which, together with derivatives of all orders, vanishes on the support of the distribution $T$. Then we have $\langle T,  R^*h \rangle=0$.  Thus, for each $h\in\mathscr{D}(X)$ with $\hbox{supp}(h)\in X\setminus\hat K$ we have $\langle RT, h\rangle=\langle T,[R^*h]\rangle=0$. This means that $\hbox{supp}(RT)\subset\hat K$.  

Before proving the second statement of  Theorem \ref{th1}, we have to show that the dual Radon transform and the convolution operation commute: 

\begin{lemma}\label{lem1}
Let $\varphi(z)\in\mathscr{D}(\mathbb{C}^n)$. Then for every $\psi(w,s)\in\mathscr{E}(X)$ the following formula holds:
$$\varphi*[R^*\psi]=R^*[\hat\varphi*_s\psi],$$
where $\hat\varphi(w,s)$ is the Radon transform of $\varphi$, and $*_s$ denotes the 
convolution with respect to the second variable $s$.   
\end{lemma}

{\it Proof.}\ For every function $\alpha(z)\in\mathscr{D}(\mathbb{C}^n)$ we have

\begin{equation}\label{g2eq7}
\int\limits_{\mathbb{C}^n}(\varphi*[R^*\psi])(z)\alpha(z)d\omega_{2n}(z)=
\int\limits_{\mathbb{C}^n} [R^*\psi](z)\left(\alpha*\varphi_1\right)(z)d\omega_{2n}(z),
\end{equation}
where $\varphi_1(z)=\varphi(-z)$. 
Let $J$ be the integral on the right-hand side of (\ref{g2eq7}). It follows from 
(\ref{g2eq2}) that 
$$J=\int\limits_{S^{2n-1}\times\mathbb{C}}\psi(w,s)\widehat{\alpha*\varphi_1}(w,s)
d\sigma(w)d\omega_2(s),$$  
where $\widehat{\alpha*\varphi_1}(w,s)$ is the Radon transform of the convolution 
$\alpha*\varphi$. We have \cite[p.p. 116-117]{Gelfand}  
$$\widehat{\alpha*\varphi_1}(w,s)=(\hat\alpha*_s\hat\varphi_1)(w,s), \quad 
\hat\varphi_1(w,s)=\hat\varphi(-w,s)=\hat\varphi(w,-s).$$
Then 
$$J=\int\limits_{S^{2n-1}\times\mathbb{C}}\psi(w,s)
(\hat\alpha*_s\hat\varphi_1)(w,s)d\sigma(w)d\omega_2(s)=
\int\limits_{S^{2n-1}\times\mathbb{C}}(\psi*_s\hat\varphi)(w,s)
\hat\alpha(w,s)d\sigma(w)d\omega_2(s).$$
In view of (\ref{g2eq2}), we have 
$$J=\int\limits_{\mathbb{C}^n}R^*[\varphi*_s\psi](z)\alpha(z)d\omega_{2n}(z).$$
Then it follows from (\ref{g2eq7}) that 
$$\int\limits_{\mathbb{C}^n}\left\{
(\varphi*[R^*\psi])(z)-R^*[\varphi*_s\psi](z)\right\}\alpha(z)d\omega_{2n}(z)=0$$
for every $\alpha(z)\in\mathscr{D}(\mathbb{C}^n)$. Therefore $(\varphi*[R^*\psi])(z)\equiv R^*[\varphi*_s\psi](z)$. The lemma is proved. 

Now suppose that the support of the Radon transform $RT$ of a distribution $T\in\mathscr{O}^{\prime}_C(\mathbb{C}^n)$ is contained in $\hat K$.  Let 
$\{\alpha_m(z)\}_{m=1}^{\infty}$ be a sequence of smooth functions on $\mathbb{C}^n$ with $\hbox{supp}(\alpha_m)\subset\{z:|z|\le1/m\}$ that converges in the space of measures to the delta function at the origin. We assume that the functions 
$\alpha_m(z)$ are even, i.e., $\alpha_m(-z)=\alpha_m(z)$. 
We set $T_m=T*\alpha_m$.  
Then the function  $T_m(z)$ belongs to $\mathscr{S}(\mathbb C^n)$ \cite[p. 244]{Schwartz}, and $T_m\to T$ in $\mathscr{O}^{\prime}_C(\mathbb{C}^n)$ \cite{Hertle1}.  Denote by $K_m$ the compact set 
$$K_m=\bigcup\limits_{z\in K}\bar B^n(z,1/m).$$
Let $\hat T_m(w,s)$ be the Radon transform of $T_m(z)$. We show that $\hbox{supp}(\hat T_m)\subset\hat K_m$. The hyperplane $\{z: \langle z, w \rangle =s\}$ meets $K_m$ if and only if there are $z^{\prime}\in K$, $z^{\prime\prime}\in\bar B^n(0,1/m)$ such that $\langle z^{\prime}, w \rangle=s-
\langle z^{\prime\prime}, w \rangle$. Therefore

\begin{equation}\label{g2eq8}
\hat K_m=\bigcup\limits_{(w,s)\in\hat K}\left(\{w\}\times\bar B^1(s,1/m)\right).
\end{equation}
Let $h(w,s)\in\mathscr{D}(S^{2n-1}\times\mathbb{C})$ be such that $\hbox{supp}(h)\cap\hat K_m=\emptyset$. Since the functions $\alpha_m$ are even, it follows from (\ref{g2eq4}) that 
$$\langle RT_m,h\rangle=\langle T_m, R^*(h)\rangle=
\langle T*\alpha_m, R^*(h)\rangle=\langle T, \alpha_m*R^*(h)\rangle.$$ 
Then by Lemma 1, we  have $\langle T, \alpha_m*R^*(h)\rangle =  
\langle T, R^*(\hat\alpha_m*_sh)\rangle$. Then

\begin{equation}\label{g2eq9}
\langle RT_m,h\rangle=\langle T, R^*(\hat\alpha_m*_sh)\rangle=
\langle RT,\hat\alpha_m*_sh\rangle.
\end{equation}
We claim that $\hat K\cap\hbox{supp}(\hat\alpha_m*_s h)=\emptyset$. Indeed,  suppose that $(w_0, s_0)\in\hat K\cap\hbox{supp}(\hat\alpha_m*_s h)$. This implies (since $\hat\alpha_m(w,s)=0$ for $|s|\ge 1/m$) that for some $s_1\in\bar B^1(0,1/m)$ we have 
$(w_0, s_0+s_1)\in\hbox{supp}(h)$. By (\ref{g2eq8}) we also have $(w_0, s_0+s_1)\in\hat K_m$,   
which contradicts that $\hbox{supp}(h)\cap\hat K_m=\emptyset$. Therefore $\hat K\cap\hbox{supp}(\hat\alpha_m*_s h)=\emptyset$, and it follows from (\ref{g2eq9}) (since $\hbox{supp}(RT)\subset\hat K$) that $\langle RT_m, h\rangle=0$. Therefore

\begin{equation}\label{g2eq10}
\hbox{supp}(RT_m)\subset\hat K_m.
\end{equation}
As remarked above, the functions $T_m(z)$ belong to $\mathscr{S}(\mathbb{C}^n)$. Then the distributions  $RT_m$ are given by the Radon transforms $\hat T_m(w,s)$ of functions $T_m(z)$.    

In view of (\ref{g2eq8}), there exist $R>0$ such that for all $m$ the sets $\hat K_m$ are contained in the set $\{(w,s): |s|\le R\}$.  Let $R_{\mathbb{R}}T_m(w,t)$ be the real Radon transform of  $T_m(z)$, that is 
$$R_{\mathbb{R}}T_m(w,t)=\int\limits_{\hbox{Re}\langle z,\bar w\rangle=t}T_m(z)d\lambda(z),$$
where $d\lambda(z)$ is the area element on the real hyperplane 
$\{z: \hbox{Re}\langle z,\bar w\rangle=t\}$. Then we have 
$$R_{\mathbb{R}}T_m(w,t)=\int\limits_{-\infty}^{\infty}\hat T_m(\bar w, t+ix)dx.$$
Since $\hat K_m\subset\{(w,s): |s|\le R\}$, it follows from (\ref{g2eq10}) that 
$R_{\mathbb{R}}T_m(w,t)=0$ for $|t|\ge R$. Then by the Helgason's support theorem, the supports of the functions $T_m(z)$ are compact.  

To complete the proof of Theorem \ref{th1}, we need the following lemma:
 
\begin{lemma}\label{lem2} Under the hypotheses and notation of Theorem \ref{th1}, there exist, for every $z_0\not\in K$, a neighborhood $V_{z_0}$ and $\delta>0$ such that the functions $T_m(z)$ vanish on $V_{z_0}$ for $m\ge1/\delta$. 
\end{lemma}

{\it Proof.}\ Fix $z_0\not\in K$. Then there exists a point $(w_0,s_0)\in S^{2n-1}
\times\mathbb{C}$ such that $\{z: \langle z,w_0\rangle=s_0\}\cap K=\emptyset$,   
$\langle z_0,w_0\rangle=s_0$ and the set $\mathbb{C}\setminus\{\langle z,w_0\rangle\}_{z\in K}$ is connected. Then $(w_0,\langle z_0,w_0\rangle)\not\in\hat K$. We set 

$$A=\left\{s\in\mathbb{C}\left|\right. (w_0,s)\in\hat K\right\},\quad
A_m=\left\{s\in\mathbb{C}\left|\right. (w_0,s)\in\hat K_m\right\}.$$
It follows from (\ref{g2eq8}) that 
$$A_m=\bigcup\limits_{s\in A}\bar B^1(s,1/m).$$
By definition of $\hat K$, for every $s\in A$ there exists $z\in K$ such that $\langle z,w_0\rangle=s$. Then $A=\{\langle z,w_0\rangle\}_{z\in K}$. Similarly $A_m=\{\langle z,w_0\rangle\}_{z\in K_m}$. Since the sets $K$ and $K_m$ are compact, it follows that the sets $A$ and $A_m$ are also compact. For some $R>0$ we have $A\cup A_m\subset\bar B^1(0,R)$. Since $\langle z_0, w_0\rangle\not\in A$, 
there is $\gamma>0$ such that $\langle z_0+\lambda, w_0\rangle\not\in A$ for every $\lambda\in\bar B^n(0,\gamma)$. Hence the convex compact set 
$\Gamma_1=\left\{\langle z,w_0\rangle, z\in \bar B^n(z_0,\gamma)\right\}$ 
and the set $A$ do not intersect. Fix $s_1\in\{s\in\mathbb{C}: |s|>R\}$. Then $s_1\in\mathbb{C}\setminus A$. Since the set $\mathbb{C}\setminus A$ is connected, there exists a broken line $\Gamma_2\subset \mathbb{C}\setminus A$ joining $s_1$ to the point $\langle z_0,w_0\rangle$. Thus  $(\Gamma_1\cup\Gamma_2)\cap A=\emptyset$. 
Then, since the sets $\Gamma_1\cup\Gamma_2$ and $A$ are compact, there exists $\delta\in(0,1)$ such that for all $m\ge 1/\delta$ we have 
$$\left\{\left(\Gamma_1\cup\Gamma_2\right)+B^1(0,\delta)\right\}\cap
\left\{A+\bar B^1\left(0,1/m\right)\right\}=\emptyset,$$
that is  $\left\{\left(\Gamma_1\cup\Gamma_2\right)+B^1(0,\delta)\right\}\cap A_m=
\emptyset$.  Put 
$$D=\left\{s\in\mathbb{C}: |s|>R\right\}\cup\left\{\left(\Gamma_1\cup\Gamma_2
\right)+B^1(0,\delta)\right\}.$$
By construction $D$ is a connected unbounded open set containing the point $\langle z_0+\lambda,w_0\rangle$ for every $\lambda\in\bar B^n(0,\gamma)$. We have by the definition of the sets $A_m$ that $(D\times\{w_0\})\cap \hat K_m=\emptyset$ for $m\ge 1/\delta$. Then it follows from (\ref{g2eq10}) that $(D\times\{w_0\})\cap\hbox{supp}(\hat T_m)=\emptyset$ for $m\ge 1/\delta$. Since the supports of $T_m$ are compact, it follows from \cite[Thm. 2]{Sekerin} that 
for every $\lambda\in \bar B^n(0,\gamma)$ and $m\ge1/\delta$ the functions $T_m(z)$ vanish on the hyperplane $\{z: \langle z,w_0\rangle=\langle z_0+\lambda,w_0\rangle\}$. Then, for every $z\in\bar B^n(z_0,\gamma)$ and $m\ge1/\delta$, we have $T_m(z)=0$. The lemma is proved. 

As mentioned above, $T_m\to T$ in $\mathscr{O}_C^{\prime}(\mathbb{C}^n)$. This means that 
\begin{equation}\label{g2eq11}
\lim\limits_{m\to\infty}\langle T_m,\varphi\rangle=\langle T,\varphi\rangle, \quad\forall
\varphi\in\mathscr{O}_C(\mathbb{C}^n),
\end{equation}
where $\mathscr{O}_C(\mathbb{C}^n)$ is the space of all infinitely differentiable functions $f$ on $\mathbb{C}^n$ for which there exist an integer $k$ such that $(1+|x|^2)^k\partial^p\bar\partial^qf(z)$ vanishes at infinity for all $p,q$ \cite[p. 173]{Horvath}. Since 
$\mathscr{D}(\mathbb{C}^n)\subset\mathscr{O}_C(\mathbb{C}^n)$,  formula (\ref{g2eq11}) holds for every $\varphi\in\mathscr{D}(\mathbb{C}^n)$. 
Let $\varphi\in\mathscr{D}(\mathbb{C}^n)$ be such that $\hbox{supp}(\varphi)\cap
K=\emptyset$. By Lemma \ref{lem2} for every $z\in\hbox{supp}\varphi$ 
there are $\delta(z)>0$ and a ball $B^n(z,\gamma(z))$ such that $T_m(z)=0$ on $B^n(z,\gamma(z))$ for $m\ge 1/\delta(z)$. Since the support of $\varphi$ is compact, it can be covered by a finite union of balls $B^n(z_k,\gamma(z_k))$, where $k=1,2\ldots, N$. Setting $\delta_0=\min\{\delta(z_k), 1\le k\le N\}$, we have $T_m(z)=0$ for $z\in\hbox{supp}(\varphi)$ and $m\ge 1/\delta_0$. Then it follows from (\ref{g2eq11}) that 
$$\langle T,\varphi\rangle=\lim\limits_{m\to\infty}\langle T_m,\varphi\rangle=0.$$
Since $\varphi\in\mathscr{D}(\mathbb{C}^n)$ is an arbitrary function such that
$\hbox{supp}(\varphi)\cap K=\emptyset$, we have $\hbox{supp}(T)\subset K$. The theorem is proved.

\end{document}